\documentclass{amsart}          
\usepackage{fullpage,eucal,amsmath,amsthm,amsfonts,verbatim,amsopn,amssymb}

\usepackage[all]{xy}

\newcommand\dmo{\DeclareMathOperator}

\renewcommand\sec{\section}
\newcommand\subs{\subsection*}

\newcommand\bea{\begin{array}{r@{,\ \ \ \ \ \ \ }l}}
\newcommand\ena{\end{array}}
\newcommand\barr{\begin{array}}
\newcommand\earr{\end{array}}
\newcommand\beq{\begin{equation}}
\newcommand\eeq{\end{equation}}
\newcommand\bqa{\begin{eqnarray*}}
\newcommand\eqa{\end{eqnarray*}}
\newcommand\bqan{\begin{eqnarray}}
\newcommand\eqan{\end{eqnarray}}
\newcommand\bit{\begin{itemize}}
\newcommand\ben{\begin{enumerate}}
\newcommand\een{\end{enumerate}}
\newcommand\eit{\end{itemize}}
\newcommand\bfi{\begin{figure}[htbp]}
\newcommand\efi{\end{figure}}
\newcommand\bce{\begin{center}}
\newcommand\ece{\end{center}}
\newcommand\bpr{\begin{proof}}
\newcommand\epr{\end{proof}}

\newcommand{\ignore}[1]{}

\newtheorem{theorem}{Theorem}[section]
\newtheorem{conjecture}[theorem]{Conjecture}
\newtheorem{lemma}[theorem]{Lemma}
\newtheorem{corollary}[theorem]{Corollary}
\newtheorem{proposition}[theorem]{Proposition}

\theoremstyle{definition}

\newtheorem{definition}[theorem]{Definition}
\newtheorem{example}[theorem]{Example}
\newtheorem{remark}[theorem]{Remark}
\newtheorem{remarks}[theorem]{Remarks}

\newtheorem{question}[theorem]{Question}

\newcommand\bthm{\begin{theorem}}
\newcommand\ethm{\end{theorem}}
\newcommand\bcn{\begin{conjecture}}
\newcommand\ecn{\end{conjecture}}
\newcommand\bla{\begin{lemma}}
\newcommand\ela{\end{lemma}}
\newcommand\bco{\begin{corollary}}
\newcommand\eco{\end{corollary}}
\newcommand\bpro{\begin{proposition}}
\newcommand\epro{\end{proposition}}
\newcommand\bdf{\begin{definition}}
\newcommand\edf{\end{definition}}
\newcommand\bex{\begin{example}}
\newcommand\eex{\end{example}}
\newcommand\brm{\begin{remark}}
\newcommand\erm{\end{remark}}
\newcommand\brms{\begin{remarks}}
\newcommand\erms{\end{remarks}}
\newcommand\bqu{\begin{question}}
\newcommand\equ{\end{question}}

\newcommand\msk{\medskip}

\newcommand\sr{\stackrel}
\newcommand\<{\langle}
\renewcommand\>{\rangle}
\newcommand\ot{\otimes}
\newcommand\op{\oplus}

\newcommand\sub{\subseteq}

\newcommand\ra{\rightarrow}

\newcommand\lra{\longrightarrow}

\newcommand\seq{\sim}
\newcommand\teq{\approx}

\newcommand\rst{|}

\newcommand\al{{\alpha}}
\newcommand\be{{\beta}}
\newcommand\ga{{\gamma}}
\newcommand\dl{\delta}
\newcommand\dltl{\tilde{\delta}}
\newcommand\lm{\lambda}
\newcommand\Lm{\Lambda}

\newcommand\omt{\omega^\tau}
\newcommand\Om{\Omega}
\newcommand\Omt{\Omega^\tau}
\newcommand\Omtp{\Omega^{\tau,p}}

\newcommand\Omtcx{\Omega^{\tau\bullet}}

\newcommand\pd{\partial}

\newcommand\xbar{\overline{x}}

\newcommand\abar{\overline{a}}

\newcommand\A{{\mathbb{A}}}

\newcommand\Z{{\mathbb{Z}}}

\renewcommand\P{{\mathbb{P}}}

\newcommand\calp{{\mathfrak{p}}}
\newcommand\calO{{\mathcal{O}}}
\renewcommand\O{{\mathcal{O}}}

\dmo{\Aff}{Aff}
\dmo{\Hom}{Hom}
\dmo{\dimn}{dim}
\dmo{\HH}{H}

\newcommand\BSpec{\textup{\textbf{Spec}}}
\dmo{\Der}{Der}
\dmo{\Ker}{Ker}
\dmo{\Mor}{Mor}
\dmo{\spec}{Spec}
\dmo{\sym}{Sym}

\newcommand\Ext{\textup{Ext}}

\newcommand\bspec{\mathop{ \mathbf{Spec} }}

\newcommand\prl{^{(1)}}

\newcommand\ME{{\mathcal{E}}}
\newcommand\MF{{\mathcal{F}}}

\newcommand\taufs{$\tau$-forms}
\newcommand\PO{PO}

\begin{document}

\title{Differentials over differential fields}
\author{Eric Rosen}
\address{Department of Mathematics\\
Massachusetts Institute of Technology\\
Cambridge, MA 02139}
\email{rosen@math.mit.edu}
\urladdr{http://math.mit.edu/~rosen/}
\thanks{This paper is an expanded version of part of my dissertation, at the University of 
Illinois at Chicago (2005).  Some of this material had been originally incorporated into 
\cite{Ros05}.}
\thanks{I am grateful to Henri Gillet, David Marker, and Thomas Scanlon for helpful discussions
on this material.}

\date{\today}

\begin{abstract}
Given an algebra $A$ over a differential field $K$, we study derivations 
on $A$ that are compatible with the derivation on $K$.  There is a
universal object, which is a twisted version of the usual module of differentials, 
and we establish some of its basic properties.  In the context of differential 
algebraic geometry, one gets a sheaf of these $\tau$-differentials which can be interpreted 
as certain natural functions on the prolongation of a variety, as studied by Buium.
This sheaf corresponds to the Kodaira-Spencer class of the variety.  
\end{abstract}

\maketitle

\sec{Introduction}

In this paper, we study derivations of algebras over differential fields 
and the associated module of differentials.  A main 
idea is to develop a theory analogous to that of K\"{a}hler differentials
for differential (commutative) algebra.  
In the usual case, given a ring $A$, an 
$A$-algebra $B$, and a $B$-module, an $A$-derivation from $B$ to $M$
is a derivation whose kernel contains $A$.  Here, we suppose that 
$A$ is in addition a differential field, and define a $\tau$-derivation 
from $B$ to $M$ to be a derivation that is compatible with the derivation 
on $A$, in a sense to be explained below.  We call the universal object 
the module of $\tau$-differentials of $B$ over $A$, and establish some
basic results about it.  

In the second part of the paper, we describe the geometric meaning of
the module of $\tau$-differentials, in the context of Buium's 
differential algebraic geometry \cite{Bui93, Bui94}, that is, algebraic geometry 
over a differential field.  Buium introduced the fundamental notion of 
the prolongation of a variety, which is a torsor of the tangent variety,
and hence an affine bundle.  Whereas a differential form on a variety $X$
can be viewed as
a regular function on the tangent variety $TX$ that is linear on each
fiber, a $\tau$-differential form on $X$ is a regular function on the 
prolongation $X\prl$ that is affine on each fiber.  And whereas 
$\bspec \sym (\Om_X)$ is $TX$, for smooth $X$, $\bspec \sym (\Omt_X)$ is a 
$(\dim X + 1)$-vector bundle over $X$, which we call the prolongation
cone of $X$, into which both $TX$ and $X\prl$ naturally embed.

Many of the results of this paper hold in the more general context of
algebras over a differential ring.  Nonetheless, we have chosen to work
over a differential field as this suffices for the applications we have
considered.  For a somewhat different, more geometric approach to some
of this material, see also~\cite{Ros05}, which was motivated by work
of Hrushovski and Itai \cite{HI} on the model theory of differential fields.

\sec{Differentials over differential fields}
\label{taualg}
In this section, we introduce and develop the theory of $\tau$-differentials,
in analogy to the usual theory of differentials (see~\cite{Eis,Mat}).  
Throughout, $(K,\dl)$ will be a differential field, and all algebras will be 
$K$-algebras.  We will also assume that $K$ contains an element
$e$ with $\dl(e) = 1$, which is necessary for some of our main results.
For example, under this assumption, for any $K$-algebra $R$, there is a 
canonical embedding of $R$ into the module $\Omt_{R/K}$ of 
$\tau$-differentials. But this is not true if the derivation on $K$ is
trivial or, more generally, if the derivation of no element is a unit.
(An example of a differential ring with this property is the fraction field of 
a polynomial ring $L[x]$, with $\dl(x)=x$ and $\dl(a) = 0$, for all $a\in L$.)

We then recall the definition of the prolongation of a $K$-algebra
(see~\cite{Joh,Bui93,Gil}) and explain some connections with 
$\tau$-differentials.
  
\bdf
Let $(K,\dl)$ be a differential ring, $R$ a $K$-algebra,  and $M$ an $R$-module.
A map $t: R \ra M$ is a {\em $\tau$-derivation (over $K$)} 
if it is a derivation and,
for all $a,b \in K$, $\dl(a)t(b)=\dl(b)t(a)$.  We often write $tr$ instead of 
$t(r)$.  

Let $\Der^\tau_K(R,M)$ denote the set of such $\tau$-derivations, 
which is an $R$-module.  
\edf

\msk

Note that any $K$-linear derivation is also a $\tau$-derivation.

\bdf
Let $R$ be a $K$-algebra.  The {\em module of $\tau$-differentials} of $R$,
denoted $\Omt_{R / K}$, is the $R$-module generated by the set 
$\{\tau(r) \mid r \in R\}$, with the relations
\ben
\item  $\tau(r+s)= \tau(r) + \tau(s)$;
\item  $\tau(rs)= r\tau(s) + s\tau(r)$;
\item  $\dl(a)\tau(b)=\dl(b)\tau(a)$;
\een
for all $r,s \in R$, and for all $a,b \in K$.  We often write $\tau r$
instead of $\tau(r)$.
The map $\tau:R \ra \Omt_{R/K}$, taking $r$ to $\tau r$ is a $\tau$-derivation,
called the {\em universal $\tau$-derivation}.  
\edf

\msk

When the differential ring $(K,\dl)$ is understood, we will usually write 
$\Omt_{R}$ instead of $\Omt_{R / K}$.

\brms
\ben
\item From the definition, for any $a,b\in K$, with $\dl(b) \neq 0$, 
one gets that $\tau a=\frac{\dl(a)}{\dl(b)}\tau b$.  
In particular, $\tau a=\tau b$ if and only if $\dl(a)=\dl(b)$ and $\tau a=0$
if and only if $a$ is a constant, i.e., $\dl(a)=0$.

For $e \in K$ such that $\dl(e)=1$, the universal derivation $\tau$ maps $K$ 
into the submodule $Q$ of $\Omt_R$ generated by $\tau e$.  
There is a natural map $\iota: R \ra \Omt_R$ taking $r$ to $r\tau e$  
mapping $R$ onto $Q$.  
This map does not depend on the choice of $e$.  
Below, Lemma~\ref{pushprf}, we give a general condition under which 
this map is injective.

\item  There is also a natural surjective map from $\lm:\Omt_R \ra \Om_{R/K}$
taking $\tau r$ to $dr$.  
Below, Lemma~\ref{pushprf} again, we show that the kernel of this map is $\iota(R)$.

\een
\erms

\bla
\label{taureps}
Given a $K$-algebra $R$, $\Omt_R$ represents the functor from
the category of $R$-modules to the category of Sets that sends 
$M \mapsto \Der^\tau_K(R,M)$.  In other words, there is a natural bijection 
$$
\Hom_{R{\textup{-mod}}}(\Omt_R,M) \equiv \Der^\tau_K(R,M).
$$
\ela
  
\bpr
As in the usual case.
\epr

\bla  
\label{tdfeq}
The $R$-module $\Omt_R$ is isomorphic to the pushout $P$ of the 
following diagram.

$$
\xymatrix{
R \otimes \Om_K \ar[r]^\alpha  \ar[d]^\dltl & \Om_R \ar[d]^\beta\\
R \ar[r] &  P
}
$$      
where $\alpha$ takes $r\ot da$ to $rda$ and $\dltl$ is the $R$-module map 
taking $r \otimes da$ to $r\dl(a)$, for all $r \in R, a \in K$.
\ela

\bpr
By definition, $P$ is isomorphic to the module $(R \oplus \Om_R) / N$,
where $N$ is the submodule of $R \oplus \Om_R$ generated by
$\{\dl(a)\oplus 0 - 0 \oplus da \mid a \in K\}$.
There is a natural surjection $\Om_R \ra \Omt_R$, taking $dr$ to $\tau r$,
whose kernel $M$ is generated by the set $\{\dl(a)db - \dl(b)da \mid  a,b \in K\}$.
Thus, to give a homomorphism $f$ from $\Omt_R$ to $P$, it suffices to give a 
homomorphism $F$ from $\Om_R$ to $P$ whose kernel contains $M$.
Let $F$ be the map that sends $dr$ to $0\oplus dr$.  
We must then show that for all $a,b \in K$, $F(\dl(a)db - \dl(b)da) = 0$ in $P$.
$$
F(\dl(a)db - \dl(b)da) = 0\op \dl(a)db - 0\op \dl(b)da
=\dl(a)\dl(b)\op 0 - \dl(b)\dl(a)\op 0 = 0
$$
Note that $f$ takes $\tau r$ to $0 \op dr$.

To prove that $f$ is an isomorphism, we construct the inverse 
$g: P \ra \Omt_R$, which we also derive from a homomorphism 
$G: R \op \Om_R \ra \Omt_R$.  Choose again $e \in K$ such that $\dl(e)=1$.
For all $r \in R$, let $G(r\op 0) = r\tau e$ and $G(0 \op dr) = \tau r$.
To show that $G$ determines a homomorphism $g: P \ra \Omt_R$,
it suffices to show that for all $a \in K$, $G(\dl(a) \op 0 - 0 \op da) = 0$.
$$
G(\dl(a) \op 0 - 0 \op da) = \dl(a)\tau e + \tau a = 
\dl(a)\tau e+\dl(e)\tau a = 0
$$
as desired.  

Note that for all $r \in R$, in $P$ we have 
$r\op 0 = r \dl(e)\op 0 = 0 \op r de$, so every element of $P$ can be written 
as a sum of elements of the form $0 \op rds$, with $r,s\in R$.  
Thus $g$ takes $0\op rds$ to $r\tau s \in \Omt_R$.
Finally, it is clear that $g$ is the inverse of $f$, so $f$ is indeed an
isomorphism.
\epr

\brm
\label{freeiota}
Identifying $\Omt_R$ with $P$ via the above isomorphism $f$,
we see that $\iota:R \ra P$ takes $r$ to $0 \op rde=r\dl(e) \op 0= r\op 0$.
Thus $\iota (R)$ is the submodule $ R \op 0\sub P$.  
Below, we give a condition under which this submodule is free.
(As far as we know, it is possible that it is always free.  This would 
be true if, for example, the assumption in Lemma~\ref{intdom} that
$R$ is an integral domain is unnecessary.)
\erm
  
Let $L$ be an extension field of a field $F$.
Recall that a set $B\sub L$ is a {\em differential basis} of $L$ over $F$
if $\{dx \mid x\in B\}$ is a basis of the $L$-vector space $\Om_{L/F}$.
In characteristic 0, a differential basis is the same thing as a 
transcendence basis (\cite{Mat}, p.\ 202).

\bla
\label{kernkern}
The kernel of the $R$-module map $\dltl:R\ot_K\Om_K\ra R$, that takes
$r\ot da$ to $r\dl(a)$, is equal to $R \ot_K M$, for $M\sub \Om_K$ 
the $K$-vector space generated by $\{\dl(a)db - \dl(b)da \mid  a,b \in K \}$. 
\ela

\bpr
We first consider the case $R=K$ is a field of finite transcendence
degree $n$. Then $\Ker(\dltl)$ is an $(n-1)$-dimensional vector space. 
It is clear that $M \sub\Ker(\dltl)$, so it suffices to show that there 
are $n-1$ linearly independent elements in $M$.
By assumption, there is an element $e \in K$ such that $\dl(e)=1$,
which must be transcendental.  Let $\{e_1, \ldots, e_n\}, e_1=e,$
be a transcendence basis for $K$, so $\{de_1, \ldots , de_n\}$
is a differential basis of $\Om_K$.  For $i=2, \ldots, n$, let
$v_i \in M$ equal $\dl(e_i)de_1 - \dl(e_1)de_i = \dl(e_i)de_1 - de_i$.
We claim that the $v_i$ are linearly independent.  
Suppose that $\sum_{i=2}^nc_iv_i=0$.  
Rearranging terms, one gets 
$
\sum_ic_iv_i = (\sum_ic_i\dl(e_i))de_1 + (\sum_i c_i de_i).
$
By the linear independence of the $de_i$, each $c_i$ must be 0, 
proving the claim.

Next, let $R=K$ be an arbitrary field.  Suppose that $\sum_ic_idm_i
\in \Ker(\dltl)$.  Let $L$ be the finitely generated subfield of $K$ 
generated by the $c_i$ and the $m_i$.  By the previous argument,
$\sum_ic_idm_i$ is contained in the $L$-vector space generated by
$\{\dl(a)db-\dl(b)da \mid a,b \in L\}$, as desired.

Finally, let $R$ be an arbitrary $K$-algebra.  By above, we have
an exact sequence of $K$-vector spaces, $0 \ra M \ra \Om_K \ra K \ra 0$.
Tensoring it with $R$, one gets the desired exact sequence of 
free $R$-modules.
\epr

\bla
\label{intdom}
Let $R$ be a $K$-algebra and an integral domain.  The map
$\al: R \ot_K \Om_K \ra \Om_R$, taking $r\ot da$ to $rda$,  
$r\in R, a \in K$, is an injection.
\ela

\bpr
Let $L$ be the fraction field of $R$.  The map $\al':R \ot_K \Om_K 
\ra \Om_L$, taking $r\ot da$ to $rda$, factors through $\alpha$,
so it suffices to show that $\al'$ is injective.

Let $B_K$ be a differential basis of $K$, and
let $B_L$ be a differential basis of $L$ such that $B_K\sub B_L$.
As $R\ot_K\Om_K$ is a free $R$-module with basis $\{db \mid b \in B_K\}$ 
and $\al'$ preserves the linear independence of the $db$, it is clear
that $\al'$ is injective, as desired.
\epr

\brm
We do not know whether the assumption that $R$ is an integral domain is
necessary.  In other words, is it true that for every $K$-algebra $R$, 
the natural map $R \ot_K\Om_K \ra \Om_R$ is injective?
\erm

\bla
\label{pushprf}  
For any $K$-algebra $R$, there is an exact sequence,
$$
R \sr{\iota}{\lra} \Omt_R \sr{\lm}{\lra} \Om_{R/K} \lra 0.
$$

Suppose that $R$ is an integral domain.  Then there is an exact sequence,
$$
0 \ra R \sr{\iota}{\lra} \Omt_R \sr{\lm}{\lra} \Om_{R/K} \lra 0.
$$
%
%
%
\ela

\bpr
The first sequence is just the pushout of the first fundamental exact
sequence along the map $\tilde{\dl}$ defined in Lemma~\ref{tdfeq}.
$$
\xymatrix{
&R \ot_K \Om_K \ar[r]^\alpha \ar[d]^{\tilde{\dl}} & \Om_R \ar[r] \ar[d]^\beta & 
\Om_{R/K}\ar[r]\ar[d]^{=} & 0 \\
&R \ar[r]^{\iota} & \Omt_R \ar[r]^{\lm} & \Om_{R/K} \ar[r] & 0
}
$$

Suppose now that $R$ is an integral domain.  By the previous lemma,
$\alpha: R\ot_K \Om_K \ra \Om_R$ is injective, 
and the second claim now follows.
\epr

\bco
\label{splitex}
Let $R$ be a smooth $K$-algebra.  Then there is a split short exact sequence,
$0 \lra R \sr{\iota}{\lra} \Omt_R \lra \Om_{R/K} \lra 0$.
\eco

\bpr
By the previous lemma and the fact that the first fundamental exact sequence 
extends to a split short exact sequence if $R$ is smooth over $K$.  
(See \cite{Mat}, p.\ 193, also for the definition of {\em smooth}).
\epr

\bco
\label{projseq}
Let $R$ be a finitely generated smooth $K$-algebra.  Then $\Omt_R$ is
a projective module.
\eco

\bpr
Since $R$ is smooth, $\Omt_R$ is locally free 
and thus projective.
Thus, by the previous corollary, $\Omt_R$ is the direct sum of two
projective modules, and thus projective.  (For the connection between
locally free and projective modules, see \cite{Eis}, Theorem A3.2.)
\epr

\msk

The following, technical lemma will be useful in the proofs, below, of the 
$\tau$-versions of the first and second fundamental exact sequences.
(See \cite{Mat}, p. 193--4.)
\bla
\label{techlem}
Let $R \ra S$ be a map between $K$-algebras.  Then $\Omt_S$ is isomorphic
to the pushout $P$ of the following diagram, where $g,h$ are the natural maps.
$$
\xymatrix{
S \ot_R \Om_R \ar[r]^g  \ar[d]^h & \Om_S\ar[d]\\
S \ot_R \Omt_R \ar[r] &  P
}
$$
\ela

\bpr
Same idea as the proof of Lemma~\ref{tdfeq}.
\epr

\msk

The next two lemmas are $\tau$-versions of basic results about
usual differentials.  They can be proved directly, as in Matsumura,
but we give different proofs, obtaining the $\tau$-sequences 
as pushouts of the usual ones.

\bla 
\textup{(First $\tau$-fundamental exact sequence)}
\label{ftfes}
Let $R \sr{f}{\ra} S$ be a map between $K$-algebras.  
There is an exact sequence of $S$-modules,
$$
S \otimes_R \Omt_R \sr{\alpha}{\lra} \Omt_S \sr{\beta}{\lra} \Om_{S / R} \lra 0
$$
where $\al(s \ot \tau r) = s\tau r$ and $\be(\tau s) = ds$.

In addition, if $S$ is smooth over $R$, then there is a short exact sequence
of $S$-modules.
$$
0 \lra S \otimes_R \Omt_R \sr{\alpha}{\lra} \Omt_S \sr{\beta}{\lra} \Om_{S / R} \lra 0
$$
\ela

\bpr
By Lemma~\ref{techlem}, we get the desired sequence as the pushout of the
first fundamental exact sequence.
$$
\xymatrix{
S \ot_R \Om_R\ar[r] \ar[d] & \Om_S \ar[r] \ar[d] & 
\Om_{S/R}\ar[r]\ar[d]^{=} & 0 \\
S\ot_R \Omt_R \ar[r]^\alpha & \Omt_S \ar[r]^\beta & \Om_{S/R} \ar[r] & 0
}
$$
The second claim follows as in Corollary~\ref{splitex}.
\epr

\bla  
\label{stfes}
\textup{(Second $\tau$-fundamental exact sequence)}  
Let $R \sr{f}{\ra} S$ be
a surjective map of $K$-algebras with $\ker(f) = I$.  Then there is an 
exact sequence of $S$-modules,
$$
I / I^2 \sr{\ga}{\lra} S \ot_R \Omt_R \sr{\al}{\lra} \Omt_S \lra 0 
$$
where $\ga(r) = 1 \ot \tau r$ and $\al(s \ot \tau r ) = s \tau r$.
\ela

\bpr
By Lemma~\ref{techlem} again, one gets the following diagram.
$$
\xymatrix{
I / I^2 \ar[r]^\mu \ar[d]^= & S \ot_R\Om_R \ar[r]^\nu \ar[d]^\eta & 
\Om_{S}\ar[r]\ar[d]^\zeta & 0 \\
I / I^2 \ar[r]^\gamma  &S \ot_R \Omt_R \ar[r]^\alpha & \Omt_S \ar[r] & 0
}
$$
\epr

\msk

The next lemma characterizes when the map $\gamma$ in the 
second $\tau$-fundamental exact sequence is injective.
(As mentioned in Remark~\ref{freeiota}, the assumption that $R,S$ are 
integral domains is perhaps unnecessary.)

\bla
Let $R \sr{f}{\ra} S$ be a surjective map of $K$-algebras, which are 
integral domains, with $\ker(f) = I$.  Then the map $\gamma$ in the 
previous diagram is injective if and only if $\mu$ is injective.
\ela

\bpr
Clearly, if $\gamma$ is injective, then so is $\mu$.
In the other direction, suppose that $\mu$ is injective, and let
$N= (I/I^2)/\Ker(\gamma)$.  Letting $M = \<\dl(a)db-\dl(b)da \mid a,b\in K \>$,
$\Ker(\zeta)=S \ot_K M$.  By the right exactness of the tensor product,
there is an exact sequence, $S\ot_R(R\ot_K M)\ra S \ot _R \Om_R
\sr{\eta}{\ra} S \Omt_R \ra 0$, so $\Ker(\eta)=J$ is a homomorphic image
of $ S\ot_R(R\ot_K M)$, which is isomorphic to $S\ot_KM$.  
By the Snake Lemma, we get the following diagram, 
with each horizontal sequence exact.
$$
\xymatrix{
0 \ar[r] & \Ker(\ga) \ar[r]\ar[d] & J \ar[r]^\zeta\ar[d] & 
S\ot_KM \ar[r]\ar[d] & 0\\
0 \ar[r] & I / I^2 \ar[r]^\mu \ar[d] & S \ot_R\Om_R \ar[r]^\nu \ar[d]^\eta & 
\Om_{S}\ar[r]\ar[d]^\zeta & 0 \\
0 \ar[r] & N\ar[r]^{\overline{\gamma}}  &S \ot_R \Omt_R \ar[r]^\alpha & \Omt_S \ar[r] & 0
}
$$
Since $J$ is a homomorphic image of $S\ot_KM$, $\zeta$ is an isomorphism,
and $\Ker(\gamma)=0$, as desired.
\epr


\msk

Given a ring $R$ and a multiplicative subset $U \sub R$, then 
$\Om_{R[U^{-1}]}=R[U^{-1}]\ot_R\Om_R$.  The next lemma establishes 
the analogous result for $\tau$-differentials.

\bla 
\textup{(Localization)}
Let $R$ be a $K$-algebra, $U$ a multiplicative subset of $R$. Then
$\Omt_{R[U^{-1}]}=R[U^{-1}]\ot_R\Omt_R$.
\ela

\bpr
By Lemma~\ref{techlem}, $\Omt_{R[U^{-1}]}$ is the pushout in the following diagram.
$$
\xymatrix{
R[U^{-1}] \ot_R \Om_R \ar[r]^g  \ar[d]^h & \Om_{R[U^{-1}]}\ar[d]\\
R[U^{-1}] \ot_R \Omt_R \ar[r]^j &  \Omt_{R[U^{-1}]}
}
$$
Since $g$ is an isomorphism, so is $j$.
\epr

\bla 
\textup{(Base Change)}
Let $R$ be a $K$-algebra, and let $(K',\dl)$ be a differential field 
extension of $K$.  Then $\Omt_{K'\ot_K R/K'} \cong K'\ot_K \Omt_{R/K}$.
\ela

\bpr
Fix $e \in K$ such that $\dl(e)=1$.
Let $T:K'\ot_KR \ra K'\ot_K\Omt_{R/K}$ be the $\tau$-derivation (over $K'$)
that maps $(a\ot r)$ to $(a\ot \tau r) + (\dl(a)\ot r\tau e)$.
By the universal property of $\tau$-differentials, this determines
a $K'\ot_K R$-module homomorphism $f: \Omt_{(K'\ot_K R)/K'}
\ra K'\ot_K\Omt_{R/K}$ that maps $\tau(a\ot r)$ to 
$(a\ot \tau r) + (\dl(a)\ot r \tau e)$.
Note also that in $\Omt_{(K'\ot_K R)/K'}$, $\tau(a\ot r)
= r\dl(a)\tau (e \ot 1) + a\tau(1 \ot r)$.

Let $U: R \ra \Omt_{(K'\ot_K R)/K}$ be the $\tau$-derivation
mapping $r$ to $\tau(1 \ot r)$.  As above, this determines 
an $R$-module homomorphism $g_0:\Omt_R \ra \Omt_{(K'\ot_KR)/K}$
mapping $\tau r$ to $\tau(1\ot r)$.  
Since $\Omt_{(K'\ot_KR)/K}$ is a $K'\ot_K R$-module, $g_0$ 
determines a $(K'\ot_K R)$-module homorphism $g:K'\ot_K \Omt_R \ra 
\Omt_{(K'\ot_KR)/K}$ taking $(a\ot \tau r)$ to $a\tau(1\ot r)$.
It is easy to see that $g$ is the inverse of $f$, so $f$ is
an isomorphism.
\epr

The following result is a $\tau$-version of a standard fact.

\bla
Let $R$ be a $K$-algebra, $f,g$ $\tau$-derivations from $R$ to $R$.
Then the commutator, $[f,g]=fg-gf$, is also a $\tau$-derivation.
\ela
  
\bpr
It is well-known that $[f,g]$ is a derivation, so it suffices to show that
for all $a,b\in K$, $\dl(a)[f,g](b)=\dl(b)[f,g](a)$.  Let $e \in K$ be such
that $\dl(e)=1$, so that for all $a \in K$, $f(a)=\dl(a)f(e)$ and $g(a)=\dl(a)g(e)$.
In fact, it is enough to show that for all $a\in K$,
$[f,g](a)=\dl(a)[f,g](e)$, since this implies that for all $a,b\in K$,
$\dl(a)[f,g](b)=\dl(a)\dl(b)[f,g](e)=\dl(b)[f,g](a)$.  
$$
\begin{array}{lll}
[f,g](a) & = & fg(a)-gf(a) \\
\\
& = & f(\dl(a)g(e)) - g(\dl(a)f(e)) \\
\\
& = & (\dl(a)fg(e)+f(\dl(a))g(e))-(\dl(a)gf(e)+g(\dl(a))f(e) )\\
\\
& = & \dl(a)[f,g](a) + \dl^2(a)f(e)g(e)-\dl^2(a)g(e)f(e)\\
\\
& = & \dl(a)[f,g](a) \\
\end{array}
$$
as desired.
\epr

\subs{Examples}

\bla
For any $n$, let $\epsilon: K[x_1, \ldots , x_n] \lra K[x_1, \ldots , x_n]$
be the map that takes any polynomial $f$ to $f^\dl$, which is obtained 
from $f$ by taking the derivative of each coefficient.  Then $\epsilon$
is a derivation on $K[x_1, \ldots , x_n]$.  Further, $\epsilon$ commutes
with each derivation $\frac{d}{dx_i}$. 

From now on, we will write $\epsilon$ as $\frac{d}{de}$, and write
$\frac{df}{de}$ for $\frac{d}{de}(f)$, which equals $f^\dl$.
\ela

\bpr
The proofs are straightforward calculations.
To simplify notation, we assume  $n=1$.  Clearly, $\epsilon$ is additive,
so it suffices to show that for any $f,g$, $\epsilon(fg)=f\epsilon(g) + 
g\epsilon(f)$.  Let $f = \sum_m a_mx^m$ and $g = \sum_n b_nx^n$.
Then
\[
\epsilon\left(\left(\sum_m a_mx^m\right)\left(\sum_n b_nx^n\right)\right) 
  = \epsilon\left(\sum_l\left(\sum_{m=0}^l a_mb_{l-m}x^l\right)\right) 
\]

\[
  = \sum_l\left(\sum_{m=0}^l(a_m\dl(b_{l-m}) + \dl(a_m)b_{l-m})\right)x^l
  =f\epsilon(g) + g\epsilon(f)
\]
as desired.

To prove that the derivations commute, it suffices to note that 
\[
\frac{d}{dx}\frac{d}{de}\sum_m a_m x^m
 = \sum_m \left(\frac{d}{dx}\dl(a_m)x^m \right) 
 = \sum_m m\dl(a_m)x^{m-1} 
 = \frac{d}{de}\frac{d}{dx}\sum_m a_mx^m.
\]
\epr

The following lemma is immediate.

\bla
\label{polyring}
Let $R = K[x_1, \ldots , x_n]$, and choose $e \in K$ such that $\dl(e) = 1$.
Then $\Omt_R$ is a free module of rank $n+1$, with generators
$\<\tau x_1, \ldots , \tau x_n, \tau e\>$.  The universal $\tau$-derivation from
$R$ to $\Omt_R$ is given by $\tau f = \sum_m \frac{\pd f}{\pd x_i} \tau x_i
+ \frac{df}{de} \tau e$.
\ela

\bla
\label{omtpoly}
Let $R = K[x_1, \ldots , x_n]$, with $e\in K$ such that $\dl (e) = 1$.
Then $\Omt_R = R \tau e \oplus (\oplus_{i=1}^n R \tau x_i)$, a free module of 
rank $n+ 1$.
\ela

\bpr
One can adapt the proof that $\Om_{R/K}$ is a rank $n$ free module, as in
\cite{Eis}.  We define $R$-module homomorphisms
$F:R^{n+1} \ra \Omt_R$ and $G: \Omt_R \ra R^{n+1}$ such that $F$ is surjective, 
and $G$ is the inverse of $F$.  Let $F$ take 
$( a_0, \ldots , a_n) \in R^{n+1}$ to 
$a_0\tau e + \sum^n_{i=1}a_i\tau x_i \in \Omt_R$.
Since  $\Omt_R$ is clearly generated
by $\{\tau e, \tau x_1, \ldots ,\tau x_n\}$, $F$ is surjective.  
  
To define $G$, note that by Lemma~\ref{taureps}
there is a natural bijection between homomorphisms from $\Omt_R$ to $R$ and
$\tau$-derivations from $R$ to itself.  Given such a $\tau$-derivation $\pd^\tau$,
let $T$ be the corresponding homomorphism.
Thus, an $(n+1)$-tuple of $\tau$-derivations determines a homomorphism 
from $\Omt_R$ to $R^{n+1}$ in an obvious way.

Let $\pd^\tau_0: R \ra R$ be the derivation from $R$ to $R$ that extends 
$\dl$ on $K$ and such that $\pd_0^\tau(x_i)=0$, for all $i$.
Clearly, $\pd_0^\tau$ is also a $\tau$-derivation.  For $m$, $1\leq m\leq n$,
let $\pd^\tau_m$ be the usual partial derivative $\pd / \pd x_m$, which is 
also a $\tau$-derivation.  Now define $G: \Omt_R \ra R^{n+1}$ by
$G(\omt) = \big(T_0(\omt), \ldots , T_n(\omt)\big)$.
It is easy to see that $G\circ F = 1_{R^{n+1}}$, the identity map on 
$R^{n+1}$, as desired.
\epr

\bpro
Let $L$ be a field extending the differential field $K$, and let $e\in K$
be an element such that $\dl(e)=1$.  Given a  set $B\sub L$,
$\{\tau b \mid b\in B\} \cup \{\tau e\}$ is a basis of the $L$-vector space
$\Omt_L$ if and only if $B$ is a transcendence basis for $L$ over $K$.
\epro

\bpr
By Lemma~\ref{pushprf}, there is an exact sequence 
$0 \ra L \sr{\iota}{\ra} \Omt_L \sr{\lm}{\ra} \Om_{L/K}\ra 0$,
with $\iota(a)=a\tau e$ and $\lm(\tau a) = da$.  Let $B\sub L$.
In one direction, suppose that $\{\tau b \mid b\in B\} \cup \{\tau e\}$ 
is a basis of $\Omt_L$.  Then $\lm(\{\tau b \mid b\in B\})=\{db \mid b\in B\}$ 
is a basis of $\Om_{L/K}$.  
Therefore, $B$ is a differential basis of $L$ over $K$,
and thus a transcendence basis.
In the other direction, if $B$ is a transcendence basis of $L$ over $K$,
then it is a differential basis, so $\lm$ maps the set
$\{\tau b \mid b\in B\} \sub \Omt_L$ bijectively onto a basis of $\Om_{L/K}$.
Thus  $\{\tau b \mid b\in B\} \cup \{\tau e\}$ is a basis of $\Omt_L$.
\epr

\bco
Let $L$ be a field extension of $K$ of transcendence degree $=n$.
Then $\dim \Omt_{L}=n+1$.
\eco

\subs{Prolongations}

Prolongations were introduced by Johnson~\cite{Joh} in the context
of what might be called differential commutative algebra.
Buium~\cite{Bui93} incorporated this work into his differential algebraic 
geometry, and developed the notion of the prolongation of 
a variety.
Here we briefly describe the algebraic version. In the next section, 
we use this to define the prolongation of a variety in the language of 
schemes.

A {\em kernel} is a ring homomorphism $f:A \ra B$ together with a 
derivation $\dl$ from $A$ into $B$.  A {\em prolongation}
is a pair of kernels, $(f,\dl):A \ra B$ and $(g,\dl ') :B\ra C$, 
such that $\dl'\circ f = g \circ \dl$.  
More generally, one can define a prolongation {\em sequence} in the obvious
way, which is how one gets, e.g., the higher prolongations of a variety.
There are also natural notions of morphisms of kernels, and of prolongations,
each of which gives a category.

The basic example of a kernel is a $K$-algebra $R$ over a 
differential field $(K,\dl)$.  A prolongation of $K\ra R$ is
a kernel $(g,\dl_S):R\ra S$ such that $\dl_S = g \circ \dl$.
A morphism between two such kernels $(g,\dl_S):R\ra S$ and
$(h,\dl_T):R\ra T$ is an $R$-algebra morphism $j:S \ra T$
such that $j\circ \dl_S = \dl_T$.  
There is a universal object in this category, which we simply
call the {\em prolongation} of $R$, and denote $R\prl$.
Given any prolongation $(g,\dl_S):R\ra S$ in the category, 
there is a unique morphism from $R\prl$ to $S$.  

\bdf
Let $R$ be a $K$-algebra.
The {\em (first) prolongation} of $R$, denoted $R\prl$, is
$\sym(\Om_R)/I$, where $\sym (\Om_R)$ denotes the symmetric algebra 
of $\Om_R$, and $I$ is the ideal generated by $\<da-\dl(a)\mid a \in K\>$.
$R\ra R\prl$ is a prolongation with the natural derivation
$\dl\prl:R \ra R\prl$, with $\dl\prl(r) = dr$, for all $r\in R$.
\edf

$\dl\prl :R\ra R\prl$ is also a $\tau$-derivation, so there is a unique 
$R$-module homomorphism $t:\Omt_R \ra R\prl$ such that 
$t\circ \tau = \dl\prl$.  Below, we show that when $R$ is smooth over $K$,
this homomorphism is an embedding.  The proof uses the following known
fact, whose geometric meaning is that the first prolongation of a smooth affine 
variety is isomorphic to the tangent variety.

\bpro
Let $R$ be a smooth $K$-algebra.  Then $R\prl \cong \sym(\Om_{R/K})$.
\epro
  
\bpr
Since $R$ is smooth, the first fundamental exact sequence, 
$0 \ra R \ot_K \Om_K \ra \Om_R \sr{\psi}\ra \Om_{R/K} \ra 0$ splits, so we can choose
a splitting homomorphism $\eta: \Om_{R/K} \ra \Om_R$.  Let $\phi:\Om_R \ra R \ot_K\Om_K$
be the map $\phi(x) = x - \eta(x)$, so we have a (non-canonical) isomorphism 
$f: \Om_R \ra (R\ot_K \Om_K) \op \Om_{R/K}$ given by $f(x) = \phi(x)\op 0 + 0\op\psi(x)$.  
This map determines an isomorphism from $\sym(\Om_R)$ to $\sym((R\ot_K\Om_K)\op \Om_{R/K})
\cong \sym(R \ot_K\Om_K)\ot_R\sym(\Om_{R/K})$, which we also call $f$.

Let $J$ be the ideal of $\sym(R\ot_K\Om_K)$ generated by the set 
$\{da - \dl(a) \mid  a \in K\}$.
The quotient $\sym(R\ot_K\Om_K) / J$ is naturally isomorphic to $R$, under the map $\dltl$
that sends each $r\in R$ to itself and sends $da$ to $\dl(a)$ for each $a \in K$.  
Tensoring the exact sequence $0 \ra J \ra \sym(R \ot_K\Om_K) \ra R \ra 0$ by $\sym(\Om_{R/K})$,
one gets the exact sequence
$$
0 \ra J\ot_R\sym(\Om_{R/K}) \ra \sym(R \ot_K\Om_K)\ot_R\sym(\Om_{R/K}) \ra R\ot_R\sym(\Om_{R/K}) \ra 0
$$
since $\sym(\Om_{R/K})$ is flat over $R$.  
Under the isomorphism  
$f: \sym(R\ot_K\Om_K)\ot_R\sym(\Om_{R/K}) \cong \sym(\Om_R)$, the ideal 
$J\ot_R\sym(\Om_{R/K})$ corresponds to the ideal $I$ of $\sym(\Om_R)$ generated by 
$\{da-\dl(a) \mid a \in K\}$, 
so $R\ot_R\sym(\Om_{R/K}) = \sym(\Om_{R/K})$ is isomorphic to $R\prl$, as desired.
\epr

\bpro
\label{embedtau}
Let $R$ be a smooth $K$-algebra. The  (unique) $R$-module homomorphism 
$t:\Omt_R \ra R\prl$ such that $t\circ\tau=\dl\prl$ is injective.
Thus, there is a canonical embedding of $\Omt_R$ into $R\prl$. 
\epro
    
\bpr
The homomorphism $t: \Omt_R \ra R\prl = \sym(\Om_R)/ I$ maps $\tau r$ to $dr + I$.
Let $t_0$ be the canonical homomorphism from $\Om_R$ to $R\prl$, taking $dr$ to $dr + I$.
Given the natural map $\beta:\Om_R \ra \Omt_R$, we have $t_0=\beta \circ t$.  Thus, to show
that $t$ is injective, it suffices to show that $\ker(t_0)=\ker(\beta)$.

Recall from the proof of Lemma~\ref{pushprf} that we have the following commutative diagram.
$$
\xymatrix{
0 \ar[r] &R \ot_K \Om_K \ar[r]^\alpha \ar[d]^{\tilde{\dl}} & \Om_R \ar[r] \ar[d]^\beta & 
\Om_{R/K}\ar[r]\ar[d]^{=} & 0 \\
0 \ar[r] &R \ar[r]^{\iota} & \Omt_R \ar[r]^{\lm} & \Om_{R/K} \ar[r] & 0
}
$$
By the Snake Lemma, $\ker(\beta) = \ker(\tilde{\dl})$, so we will show that
$\ker(t_0)=\ker(\dltl)$.  (By Lemma~\ref{kernkern}, $\ker(\dltl)$ is $R\ot_KM$, 
$M\sub\Om_K$ the $K$-vector space generated by $\{\dl(a)db - \dl(b)da\mid  a,b \in K \}$,
though we do not use this here.)

We now calculate $\ker(t_0)$.  Let $h:R\prl \ra \sym(\Om_{R/K})$ be the isomorphism 
from the previous proposition, and define $\hat{t}_0=h\circ t_0$, so 
$\ker(\hat{t}_0) = \ker(t_0)$.  The map $\hat{t}_0$ is the composite of the maps
$$
\Om_R \ra \sym(\Om_R)\sr{f}{\ra}\sym(R\ot_K\Om_K)\ot_R\sym(\Om_{R/K})\ra 
R \ot_R\sym(\Om_{R/K})\ra \sym(\Om_{R/K}).
$$
Explicitly, $m \in \Om_R$ is sent to $(\dltl(\phi(m))\ot 1) + (1 \ot \psi(m)) 
\in R\ot_R \sym(\Om_{R/K})$, with $\dltl(\phi(m))\in R$ and $\psi(m)\in \Om_{R/K}$.
Thus, $\hat{t}_0(m)=0$ if and only if $\psi(m)=0$ and $\dltl(\phi(m))=0$.
Here, $\psi(m)=0$ if and only if $m\in R\ot_K\Om_K$ and, in case $\psi(m)=0$, then $\phi(m) = m$.
Thus, $m \in \ker(\hat{t})$ if and only if $m \in R\ot_K\Om_K$ and $m \in \ker(\dltl_0)$.  
This completes the proof.
\epr

\sec{Varieties, prolongations, and $\tau$-differential forms}
\label{scheme}

In this section, we introduce the sheaf of $\tau$-differential forms
on a variety over a differential field, and describe the connection to 
the {\em prolongation} of the variety, introduced by Buium (see ~\cite{Bui93}).
First we describe the construction of the prolongation, which is a 
torsor under the tangent bundle, and thus an affine bundle.

We adopt the following conventions. $(K,\dl)$ is an algebraically closed 
differential field with an element $e \in K$ such that $\dl(e)=1$.  A 
variety is an integral, separated $K$-scheme of finite type.  We will only 
consider smooth, i.e., nonsingular, varieties (see~\cite{Har}, p.\ 268).

\subs{Affine bundles}
\label{schbun}

Recall that an affine space is a principal homogeneous space of (the additive
group of) a vector space.  In other words, given a field $K$, a $K$-affine space is
a triple $(A,V,\al)$, where $A$ is a set, $V$ a $K$-vector space, and $\alpha$
a regular action of $V$ on $A$, though we generally omit explicit 
mention of the function $\al$.
We will say that the dimension of $A$ is just the dimension of $V$.  

An affine map between $K$-affine spaces $(A,V)$ and $(B,W)$ is a
function $f: A \ra B$ such that there is a linear map $\lm f: V \ra W$
such that for all 
$a \in A, v \in V$, $f(v\cdot a) = (\lm f (v))\cdot f(e)$.
There is also a natural `linearization' functor $\lm$ from the category of
affine spaces to vector spaces, with $\lm(A,V) = V$ and behaving on 
morphisms as described above.

Given a $K$-affine space $(A,V)$, there is an associated `dual' {\em vector} 
space $(A,V)^\vee$ of affine maps from $A$ to $K$, of dimension $\dim A + 1$.

An affine bundle over a variety can then be defined in analogy to the definition
of a vector bundle (e.g., see \cite{Har}, p. 128).
           
\bdf  
Let $Y$ be a variety. A {\em (geometric) affine bundle} of rank $n$ over $Y$
is a variety $X$ with a morphism $f: X\lra Y$, together with the data of an open
covering $\{U_i\}_{i\in I}$ of 
$Y$ and isomorphisms $\psi_i:f^{-1}(U_i)\lra \A^n_{U_i}$
such that for all $i,j\in I$ and open affine $V=\spec B \sub U_i \cap U_j$, the 
automorphism $\psi=\psi_j\circ\psi_i^{-1}$ of $\A^n_V=\spec B[x_1,\ldots, x_n]$ 
is given by an affine automorphism $\theta$ of $B[x_1, \ldots, x_n]$, i.e.,
$\theta(b)=b$ for $b \in B$ and $\theta(x_j)=c_j+\sum_kb_{jk}x_k$,
for $c_j \in B$ and suitable $b_{jk} \in B$.    
\edf   
    
\brm
Given a rank $n$ affine bundle over a variety $Y$, and a point 
$\calp\in\spec B \sub Y$, the fiber $Y_\calp$ has the structure of an
$n$-dimensional affine space over $\kappa(\calp)=B_\calp / \calp_\calp$.
\erm

\bdf     
A morphism of affine bundles, $f:X\lra Y$ and $g:W\lra Z$, is a pair of morphisms
$s:X\lra W, t:Y\lra Z$ so that $t \circ f = g \circ s$ and, for any $a \in Y,$ 
$b = t(a) \in Z$, there are affine neighborhoods $U = \spec A$ of $a$, 
$V = \spec B$ of $b$, so that
$$
f^{-1}(U) \cong \spec A [ x_1, \ldots, x_n]
$$
$$
g^{-1}(U) \cong \spec B [ y_1, \ldots, y_m] 
$$
and the induced map $s:f^{-1}(U)\lra g^{-1}(V)$ is given by an affine
map $h:B[y_1,\ldots, y_m]\lra A[x_1,\ldots ,x_n]$ such that 
(i) $h$   maps $B$ to $A$ and $h|_B = t$ and (ii) $h(y_i)=a_i + \sum_j c_{ij}x_j$.
\edf  

Equivalently, one can define an affine bundle as a {\em torsor} of a vector bundle.

\ignore{
Recall the definition of a torsor.  Let $Y,Z$ be varieties such that 
$Z \lra Y$ is a vector bundle (or more generally a group scheme).  
A torsor under 
$Z\lra Y$ is a morphism $X \lra Y$ and an action $Z\times_Y X \lra X$
such that the induced morphism $Z\times_YX\lra X\times_Y X$ is an 
isomorphism.  A torsor is {\em trivial} if it is isomorphic to the group scheme
itself (if and only if it has a section) and {\em locally trivial} if $Y$ can be covered
by Zariski open sets over which the torsor is trivial.
%
%
The fibers of a torsor under a vector bundle clearly have the structure
of affine spaces, in the sense of Definition~\ref{aff}.  When in addition
the torsor is locally trivial, then it is an {\em affine bundle}.
(Compare the definition of a vector bundle in~\cite{Har}, p.\ 128.)
}

\subs{Kernels and prolongations}
\label{schprol}

We now describe Buium's `globalization to the frame of schemes'
of Johnson's work.  

\bdf
Let $X$ be a scheme, $\MF$ a sheaf of modules on $X$.  A derivation $\dl$
from $\O_X$ to $\MF$, is a set of derivations $\dl_U:\O_X\rst_U\ra \MF\rst_U$,
for $U \sub X$ open, compatible with the restriction maps.
Let $\Der(\O_X,\MF)$ denote the set of such derivations.
Likewise, let $\Der^\tau(\O_X,\MF)$ denote the set of $\tau$-derivations
from $\O_X$ to $\MF$, defined in the obvious way.

Given a morphism of schemes $g:Y \ra X$, a derivation from $X$ to $Y$,
also written $\dl: X \ra Y$, is a derivation $\dl \in \Der(\O_X,g_*\O_Y)$.
Similarly for $\tau$-derivations.
\edf

\brm
\label{basicex}
For any $K$-variety $X$, the basic example of a ($K$-linear) derivation is 
the differential map $d: \O_X \ra \O_{TX}$.  The map $d$ determines an embedding
of $\Om_{X/K}$ into $p_*\O_{TX}$, $p:TX \ra X$, as $\O_X$-modules.
One can also consider $\Om_{X/K}$ as a sheaf of abelian groups on $TX$,
which is not, however, an $\O_{TX}$-module.

Below, we will see that the map $\dl\prl$  from $\O_X$ to $\O_{X\prl}$, 
is a $\tau$-derivation, and closely related to $d$.  In particular, over
a field $K$ with a trivial derivation, $X\prl = TX$, and $\dl\prl=d$.
\erm

The following lemma is clear.

\bla
\label{natbij}
Let $R$ be a $K$-algebra and $M$ an $R$-module, so $X = \spec R$ is 
an affine scheme and $\MF = (M)^\sim$ is an $\O_X$-module.  There are
natural bijections between $\Der(R,M)$ and $\Der(\O_X,\MF)$, and also
between $\Der^\tau(R,M)$ and $\Der^\tau(\O_X,\MF)$.
\ela

\bla
\label{taulocal}
Let $X$ be a variety, $\MF$ a quasi-coherent sheaf on $X$, and $t$ a map 
from $\O_X$ to $\MF$.  To prove that $t$ is a ($\tau$-) derivation, it
suffices to check that for any affine $Y \sub X$, $Y = \spec R$,
$\MF\rst_Y = M^\sim$, $M$ an $R$-module, the map
$t\rst_Y:R \ra M$ is a ($\tau$-) derivation.
\ela

\bpr
Straightforward.
\epr

\msk

We can now give Buium's definitions of kernel and prolongation for varieties.
For {\em affine} varieties, these are obviously equivalent to Johnson's.

\bdf
A {\em kernel} from a variety $Z$ to a variety $X$ is a pair $(f,\dl)$,
$f$ a morphism from $Z$ to $X$ and $\dl \in \Der(\O_X,f_*\O_Z)$.
\edf

There is also a `relative' notion of a kernel for $K$-varieties, where
$\dl$ must be compatible with the derivation on $K$.  This makes $\dl$
into a $\tau$-derivation.
For our purposes, the relative version is more important, and the only one
that we will consider.  Given any $K$-variety $X$, there is a natural
kernel $X \ra K$, as well as the fundamental example of the kernel from 
$X\prl$ to $X$, defined below, which prolongs the former.

We first define the prolongation of an affine variety, and then show 
how to globalize.

\bdf
\label{affprol}
Let $X$ be an affine $K$-variety, $X = \spec R$, $R$ a $K$-algebra.
The first prolongation $X^{(1)}$ is $\spec \sym (\Om_R) / I$,
where $I$ is the ideal generated by $\{da - \dl(a) \mid  a \in K\}$.

The projection from $X^{(1)}$ to $X$ is determined by the 
natural embedding of $R$ into $\sym (\Om_R) / I$.  By Lemma~\ref{natbij},
the $\tau$-derivation from $\dl\prl: R \ra R\prl$ corresponds to 
a $\tau$-derivation $\dl\prl:X\ra X\prl$, making $X\prl\ra X$ into a kernel.
\edf

To see that $X^{(1)}$ is a torsor under the tangent bundle, 
it is useful to give an equivalent definition, in terms of representable
functors (see, e.g., \cite{EH}, sections I.4 and VI.1).  
This approach also globalizes more easily, that is, without explicit patching,
and provides some insight into the connection between prolongations 
and $\tau$-derivations.

\bdf
\label{affproldef}
Let $X$ be an affine $K$-variety, $X = \spec R$, $R$ a $K$-algebra.
By smoothness, there is an exact sequence of $R$-modules,
$0\lra R\otimes_K \Om_K \lra \Om_R \lra \Om_{R/K} \lra 0$.
Let $\dltl: R\otimes_K \Om_K \lra R$ be the $R$-module homomorphism given by
$\dltl(r\otimes da) = r\cdot \dl(a)$.  Let $F$ be the functor from the category
of $R$-algebras to the category of Sets that associates to any $K$-algebra $S$, 
the set of pairs $(g,w)$, $g:R\ra S$ a $K$-algebra map and $w:\Om_R \ra S$
a $K$-module homomorphism making the following diagram commutate.
$$
\xymatrix{
0 \ar[r] & R\otimes_K \Om_K \ar[r]^\alpha \ar[d]^\dltl & \Om_R \ar[r]^\beta \ar[d]^w& 
\Om_{R/K} \ar[r] & 0\\
& R \ar[r]^g & S  & &
}
$$
Let $R^{(1)}$ be the $K$-algebra representing this functor, and let
$X^{(1)}=\spec R^{(1)}$.  
\edf

\brm
It is easy to check that $R^{(1)}$ is isomorphic to $\sym (\Om_R) / I$,
as in Definition~\ref{affprol}.
\erm

We now give the general definition of the prolongation of a $K$-variety $X$
and show that it is a $TX$-torsor.  (See~\cite{Bui93}, p.\ 1392--93.)

\bdf  
\label{prolfin}
Let $X$ be a $K$-variety, $f:X \ra K$ the structure map.  
There is an exact sequence
$$
0 \lra f^*\Om_K \sr{\alpha}{\lra} \Om_X \sr{\beta}{\lra} \Om_{X/K} \lra  0
$$
of $\O_X$-modules.  Let $\dltl: f^*\Om_K \ra \O_X$ be the map determined
by $\dltl$ from Definition~\ref{affproldef}.  Let $F$ be the functor
from $K$-schemes to Sets that takes a scheme $Z$ to the set of pairs
$(g,w)$, $g:Z \ra X$ a morphism of schemes, $w:\Om_X \ra g_*\O_X$ a map
of $\O_X$-modules such that the following diagram commutes.
$$
\xymatrix{
f^*\Om_K \ar[r]^\alpha \ar[d]^\dltl & \Om_X \ar[d]^w \\
\O_X \ar[r]^g & g_*\Om_Z
}
$$
The prolongation of $X$, written $X\prl$, is the scheme representing $F$.

Explicitly, $X\prl$ is the scheme $\BSpec\,  \sym(\Om_X)/I$, where $I$ is the 
ideal sheaf in the symmetric algebra $\sym(\Om_X)$ generated by all
elements of the form $\dltl(x)-x$, $x$ a local section of $f^*\Om_K$
(as in Definition~\ref{affprol}).

Let $g\prl:X\prl\ra X$ be the natural morphism of schemes, and $w\prl:\Om_X\ra
(g\prl)_*\O_{X\prl}$ the natural map of $\O_X$-modules.  $w\prl$ determines
a $\tau$-derivation in $\Der^\tau(\O_X,g\prl_*\O_{X\prl})$, 
which we call $\tau\prl$, though Buium, who introduced it, called it 
$\tilde{\delta}$ (\cite{Bui93}, p. 1396).
This makes $X\prl$ into a kernel 
$(g\prl,\tau\prl):X\prl \ra X$, which prolongs the kernel $X\ra K$.
\edf

\bdf
With the notation of the previous definition, let
$G$ be the functor from $K$-schemes to Sets that takes a scheme
$Z$ to the set of pairs $(g,v)$, $g:Z \ra X$ a morphism of schemes
and $v: \Om_{X/K} \ra g_*\O_Z$ a map of $\O_X$-modules.
$G$ is represented by the tangent variety $TX$ (which can also be 
constructed as $\BSpec\, \sym(\Om_{X/K})$).
  
Considering schemes as the functors they represent,
we have functorial transformations $X\prl \ra X$ and $TX \ra X$,
and an action $TX \times_X X\prl \ra X\prl$ given by
$$
((g,v),(g,w))\mapsto (g, w + v\circ \beta), (g,v)\in G(Z),(g,w)\in F(Z)
$$
This makes $X\prl$ into a $TX$-torsor.
\edf

Let $X,Y$ be varieties, and $f:X\lra Y$ a morphism between them.
We recall how the lifting map $f^{(1)}:X^{(1)}\lra Y^{(1)}$ is defined,
and that it is compatible with the torsor structure.
(Compare~\cite{Bui93}, p.\ 1435, where Buium writes that $f\prl:X\prl\lra Y\prl$ is
`equivariant' with respect to $df:TX \lra TY$.)
Given, $f:X\lra Y$ the lifting $f^{(1)}: X\prl \lra Y\prl$ determines
a morphism in the category of kernels.  
In fact, given $f$, $f^{(1)}$ is the unique morphism 
from $X^{(1)}$ to $Y^{(1)}$ so that $f,f^{(1)}$ is a morphism of kernels.
It suffices to consider the affine case, so assume $X =\spec T,
Y = \spec R$, and that the morphism $f:X \lra Y$ corresponds to a ring
homomorphism $R \lra T$, which we also denote by $f$.

We have $R^{(1)} = \sym(\Om_R)/I_R$ and $T^{(1)} = \sym(\Om_T)/I_T$.  Then
$f^{(1)}: R^{(1)} \lra T^{(1)}$ is given by, for $b \in R$,
$f^{(1)}(b) = f(b)$ and $f^{(1)}(db) = d(f(b))$.
This yields the following commutative diagram.
$$
\xymatrix{
T\prl & R\prl \ar[l]_{f\prl}\\
T\ar[u]^{(p_T,\dl_T)}&R\ar[l]_{f}\ar[u]_{(p_R,\dl_R)}
}
$$
Note that we also get that $\dl_T\circ f = f\prl\circ\dl_R$, which is
precisely the condition for having a morphism of kernels.

Next one wants to show that $f\prl$ is compatible with the torsor structure,
that is, that the following diagram is commutative.
$$
\xymatrix{
X\prl\times_X TX\ar[r]^{f\prl\times Tf} \ar[d]_{m_X} & Y\prl\times_YTY
\ar[d]^{m_Y}\\
X\prl\ar[r]^{f\prl} & Y\prl
}
$$
This diagram corresponds to:
$$
\xymatrix{
\sym(\Om_T/I_T)\otimes\sym(\Om_{T/K}) & \sym(\Om_R/I_R)\otimes\sym(\Om_{R/K})
\ar[l]_{f\prl\otimes df}\\
\sym(\Om_T)/I_T \ar[u]^{m_T} & \sym(\Om_R)/I_R\ar[u]^{m_R}\ar[l]_{f\prl}
}
$$
It now suffices to observe that the following diagram is commutative.
$$
\xymatrix{df(a)\otimes 1&da\otimes 1\ar[l]\\
df(a)\ar[u]&da\ar[l]\ar[u]
}
$$

\subs{$\tau$-differentials on schemes}
\label{tausch}

\ignore{
In Section~\ref{tauforms}, given a variety $X$, we defined a map 
$\tau: K[X] \ra \Omt[X]$ and showed, in Lemma~\ref{twder}, that $\tau$ is 
a derivation.  Recall that the map $\tau$ is a twisted version of the 
map $d: K[X] \ra K[TX]$ that takes a regular function $f$ on $X$ to its
differential $df$, considered as a regular function on $TX$.
Here, we develop this material scheme-theoretically, and provide some 
additional information.
}

The universal derivation from a $K$-algebra $R$ to $\Om_{R/K}$ corresponds,
geometrically, for a $K$-variety $X$, to a derivation in $\Der(\O_X, \Om_{X/K})$.
Alternatively, it can be considered as belonging to $\Der(\O_X, p_*\O_{TX})$,
$p:TX \ra X$, but this is really equivalent, as there is a canonical 
$\O_X$-module embedding of $\Om_{X/K}$ in $p_*\O_{TX}$.  
Thus $\Om_{X/K}$ is also naturally a sheaf of abelian groups on $TX$.
In differential algebraic geometry, there is twisted version of this 
picture, described below.

Let $f: X \ra Y$ be a morphism of varieties, and $f^{(1)}:
X^{(1)} \ra Y^{(1)}$ the lifting morphism on their prolongations.
For $Y = \A$, the affine line, one can modify $f^{(1)}$ to give
element of $\O_{X^{(1)}}(X\prl)$, as we now describe.
Thus we will get a map 
$\tau:\calO_X(X) \calO_{X\prl} (X\prl)$.
The natural bijection between $\Mor(X,\A)$ and $\O_X(X)$ takes $(f,f^\#)$, 
$f^\#:K[x]\ra\O_X(X)$ a ring homomorphism, to $f^\#(x)\in\O_X(X)$.
Let $\tau_0:\calO_X(X) \ra \Mor(\calO_{X\prl},\A\prl)$ be the map 
$f \mapsto f\prl$.  We have $\A\prl\cong T\A=\spec K[x,\dl x]$, so let
$q:\A\prl\ra\A$ be the projection onto the fiber ({\em not} the base), given by the 
ring embedding $K[\dl x]\ra K[x,\dl x]$.  Thus $q$ induces a map
from $\Mor (\calO_{X\prl},\A\prl)$ to $\Mor(\calO_{X\prl},\A)
\cong\calO_{X\prl}(X\prl)$.  Putting the pieces together, 
$\tau := q \circ \tau_0$ is the desired map from $\calO_X(X)$ to
$\calO_{X\prl} (X\prl)$.
  
One can also consider rational maps on $X$, that is, morphisms from open 
subschemes of $X$ to $\A$.  For each such $Y \sub X$, one gets a map 
$\tau_Y$ from $\O_Y(Y)$ to $\calO_{Y\prl} (Y\prl)$, as above.
Equivalently, the $\tau_Y$'s determine a map, which we also call $\tau$,
between the $\O_X$-modules $\O_X$ and $p_*\O_{X\prl}$, for $p$ the 
canonical projection from $X\prl$ to $X$.  

\bla
Given a variety $X$, the map $\tau: \O_X \ra p_*\O_{X\prl}$
is a $\tau$-derivation, that is, 
$\tau \in \Der^\tau(\O_X,p_*\O_{X\prl})$.
\ela

\bpr
By Lemma~\ref{taulocal}, it suffices to prove that for all 
affine subschemes $Y \sub X$, $Y \cong \spec T$, then
$\tau_Y: T \ra T\prl$ is a $\tau$-derivation.

Passing to the category of $K$-algebras, and letting $R = K[x]$,
$R\prl=K[x,\dl x]$ one gets the following diagram.
$$  
\xymatrix{
T\prl & R\prl \ar[l]_{f\prl} & K[\dl x] \ar[l]_{\pi} \\
T \ar[u]^{(p_T,\dl_T)} & R \ar[l]_{f} \ar[u]_{(p_R,\dl_R)}
}
$$
Given $f \in T=\calO_X(X)$, $\tau(f)=f\prl\circ\pi(\dl x)$, so we get
$\tau(f)=\dl_T(f)=df\in\sym(\Om_T)/I_T$, which is easily seen to be a  
$\tau$-derivation.
\epr

\bpro
\label{coincide}
Let $X$ be a variety.  The $\tau$-differential map $\tau:\O_X\ra p_*\O_{X\prl}$
is the same as Buium's map, $\tau\prl:\O_X\ra p_*\O_{X\prl}$.
\epro

\bpr
Immediate from the proof of the previous lemma and Definition~\ref{prolfin}.
\epr

\msk

We now define the coherent sheaf of $\tau$-differentials on a variety $X$,
written $\Omt_X$, which will be locally free of rank dim$(X)+1$.
As a subsheaf of $p_*\O_{X\prl}$, these can be viewed as rational 
functions on $X\prl$, which are affine maps on each fiber of $X\prl\lra X$.

For $X = \spec T$ affine, let $\hat{\Om}^\tau_T$ be the submodule
of $T\prl$ generated by $\dl_T(T)$ (which is isomorphic to $\Omt_T$,
by Proposition~\ref{embedtau}).  Then $\Omt_X$ is the $\O_X$-sheaf 
$(\hat{\Om}^\tau_T)^\sim$, which naturally embeds in the
$\O_X$-sheaf $(T\prl)^\sim$.  The following result globalizes this 
to varieties.

\bla
\label{projtauseq}
Let $X$ be a $K$-variety, $\pi:X \ra K$.
There is the following diagram of sheaves on $X$, with each row exact.
In particular, $\Omt_X$ is locally free of rank $\dim X + 1$.
$$
\xymatrix{
0 \ar[r] &  \pi^*\Om_K \ar[d]^{\dltl} \ar[r] & \Om_X \ar[r] \ar[d] 
& \Om_{X/K}
\ar[r] \ar[d]^= & 0 \\
0 \ar[r] & \calO_X \ar[r] & \Om^\tau_X \ar[r]^{\Lm} & \Om_{X/K} \ar[r]  & 0
}
$$    
\ela

\bpr
From the proof of Lemma~\ref{pushprf} and Corollary~\ref{splitex}.
\epr

As an extension of $\Om_X$ by $\O_X$, $\Omt_X$ corresponds to an element
of the cohomology group $\Ext^1(\Om_X,\O_X)$, which is naturally isomorphic
to $H^1(X,\Theta_{X/K})$, the dual sheaf of $\Om_{X/K}$.  In fact, $\Omt_X$
corresponds to the Kodaira-Spencer class of $X$, as defined in \cite{Bui93},
p.\ 1396, as can be easily seen by comparing the diagram in the preceding lemma
with Buium's construction.  (See also \cite{Ros05}).

\bpro
Given a variety $X$, the sheaf $\Omt_X$ corresponds to the Kodaira-Spencer
class of $X$.
\epro

\bla
Let $f:X \ra Y$ be a morphism of $K$-varieties.  Then there is an
exact sequence of sheaves on $X$,
$$
f^*\Omt_Y \lra \Omt_X \lra \Om_{X/Y} \lra 0
$$
\ela

\bpr
From Lemma~\ref{ftfes}.
\epr

\brm
Everything is functorial so, for example, given a morphism of $K$-varieties,
$f:X\lra Y$, the following diagram of $\O_X$-modules commutes.  
$$
\xymatrix{
\Om_X^\tau \ar[d] & f^*\Om^\tau_Y \ar[l]\ar[d] \\
\Om_{X/K} & f^*\Om_{Y/K}\ar[l]
}
$$
\erm

\section{The prolongation cone}
\label{cone}

We introduce a new construction, the prolongation cone of a variety.  
If the variety $X$ is smooth, then it will be the smallest vector bundle over $X$ 
into which both $TX$ and $X\prl$ can be embedded.  

\bdf
Let $X$ be a $K$-variety and $\Omt_X$ the sheaf of $\tau$-differentials.
The {\em prolongation cone} of $X$, written $CX$, is 
$\bspec \sym (\Omt_X)$.
\edf

For $X$ smooth, $\Omt_X$ is a locally free sheaf, and $CX$ is the 
geometric vector bundle associated to it.

To prove that there is a closed embedding of $X\prl$ in $CX$ over $X$, 
it suffices to prove this on affine subvarieties of $X$.  
Reformulated in terms of $K$-algebras, this is equivalent to showing 
that, for any $K$-algebra $R$, there is a natural surjective $R$-algebra 
homomorphism from $\sym(\Omt_R)$ onto $R\prl$.

\bpro
Let $R$ be a $K$-algebra, and $\tilde{R} := \sym(\Omt_R)$.
There is a natural surjective $R$-algebra homomorphism $\tilde{R} \ra R\prl$.
\epro

\bpr
Let $S$ be the $R$-algebra $S:= \sym (\Om_R)$.  Recall that the prolongation
of $R$ is $R\prl : = S / J$, $J$ the ideal generated by 
$\< da - \dl a \mid a \in K \>$.  Let $\tilde{R} := \sym (\Omt_R)$.

Let $f:S \ra R\prl$ be the natural quotient map, with kernel $J$.
Since there is a natural surjection of $R$-modules from 
$\Om_R$ to $\Omt_R$, there will also be a natural surjection $g$ from
$S$ onto $\tilde{R}$, with some kernel $I$.  Thus, to give the desired
surjection from $\tilde{R}$ onto $R\prl$, it suffices to show that
$I \sub J$.  Then $f$ will factor through $g$, so that the desired
surjection $h:\tilde{R} \ra R\prl$, $f = h \circ g$, has kernel $\cong J/I$.

The kernel $N$ of the natural $R$-module homomorphism from 
$\Om_R$ to $\Omt_R$ is generated by $\<\dl(a)db-\dl(b)da \mid  a,b \in K\>$
(which is immediate from the definition of $\Omt_R$,
 but see also Lemma 5.7 and the proof of 
Proposition 5.25).  Thus the kernel $I$ of $g:S\ra\tilde{R}$ is the 
ideal $I : = 0 \op N \op N^2 \op \ldots$.

Finally, we observe that $I \sub J$.  The ideal $I$ is generated
by elements of the form $\dl(a)db - \dl(b)da, a,b \in K$, which can be 
written as $da(db-\dl(b))-db(da-\dl(a))$, and are thus in $J$.
\epr

\bco
For any variety $X$, there is a natural closed embedding of 
$X\prl$ into $CX$.
\eco

For any $K$-variety $X$, let $TX$ denote $\bspec \sym (\Om_{X/K})$,
which equals the usual tangent variety of $X$ when $X$ is smooth.
To prove that there is a closed embedding of $TX$ in $CX$, one can
again, as above, reduce it to an assertion about $K$-algebras.

\bpro
Let $R$ be a $K$-algebra.  Then there is an $R$-algebra homomorphism
from $\sym (\Omt_{R})$ onto $\sym (\Om_{R/K})$.
\epro

\bpr
The $R$-module homomorphism from $\Omt_R$ onto $\Om_{R/K}$ determines
such a map.  Alternatively, one can argue as in the previous proposition.
\epr

\bco
For any variety $X$, there is a natural closed embedding of 
$TX$ into $CX$.
\eco

Let $X$ be an affine variety, together with a closed embedding $X \ra \A^n$.
We now describe $CX$ and the embeddings of $X\prl$ and $TX$ in 
$CX$ in local coordinates.

Above we defined an affine space as a principal homogeneous space of a 
vector space.
For the next proposition, it will be helpful to recall an alternative
characterization of an affine space as a coset $A$ of a linear
subspace of a vector space $V$.  In this case, say that $A$ is a {\em proper}
affine space if $0 \not\in A$.  Given a proper $A$, let 
$A^\circ := \{ca \mid  c \in K, a \in A\}$, the smallest linear space
containing $A$.  If $A$ is proper, there is a surjective homomorphism
$V^\vee \ra \Aff (A)$.  If, additionally, $A$ has codimension $1$, then this
is an isomorphism.  
Note that, given an affine space $A \sub V$, there is an associated, 
isomorphic, affine space $A'\sub V' := V \times K$, $A' := \{a \times 1 \mid
a \in A\}$, that is clearly proper.  In particular, if $A = V$, then
$A' \sub V'$ has codimension $1$ and $V'^\vee \cong \Aff(A)$.

In the same way, one can also define an affine bundle over a variety $X$ to be
a closed subvariety of a vector bundle $Y$ over $X$ with the obvious properties.

\bpro
Let $X = \spec B \sub \A^n$ be an affine variety.  Let $A:=K[x_1,\ldots, x_n]$
and $B = A / I$, $I$ an ideal of $A$.  Let $C = K[x_1,\ldots x_n, \tau x_1,
\ldots , \tau x_n, \tau e]$.  Then $CX = \spec C/J$,
where $J$ is the ideal generated by $I$ and $\{\tau f \mid f \in I \}$.
\epro

\bpr
First, we calculate $X\prl$, for $X = \A^n$.  Of course, $X\prl = TX =
\A^{2n}$, but we want to consider $X\prl$ as an affine bundle, so we embed
$X\prl$ in $\A^{2n+1} : = \spec K[x_1, \ldots , x_n, \tau x_1, \ldots ,
\tau x_n, \tau e]$ by sending $\abar \in \A^{2n}$ to 
$\abar \times 1 \in \A^{2n+1}$.  That is, for each $b \in X$, $X\prl_b$ is the 
associated affine space of the vector space $TX_b \sub \A^n$.

Thus, affine bundle maps on $X\prl$ can be written as 
$\sum_if_i(\xbar)\tau x_i + g(\xbar) \tau e$.  
(See also Lemma~\ref{polyring}.)
In particular, for each $f \in A$,
$\tau f$ is an affine bundle map on $X\prl$. For each $K$-valued 
point $b \in X$, $\tau f(b)$ is an affine map on the affine space
$X\prl_{b}$.

We now consider the general case, $X = \spec B, B = A/ I$.
By Lemma~\ref{stfes}, we have an exact sequence of $B$-modules,
\[
\xymatrix{
I / I^2 \ar[r]^\alpha & B \ot_A\Omt_A \ar[r]^\beta & \Omt_B \ar[r] & 0
}
\]
where $\alpha(f) = 1 \ot \tau f$, for $f \in I$.  The map $\alpha$ is not
necessarily injective, but there is always the related short exact sequence,
\[
\xymatrix{
0 \ar[r]& \alpha(I / I^2) \ar[r]& B \ot_A\Omt_A \ar[r]^\beta & \Omt_B \ar[r] & 0.
}
\]
where $\alpha (I / I^2) = \{1 \ot\tau f \mid f \in I\} \sub B \ot_A\Omt_A $.

Since $\Omt_A$ is a free rank $n+1$ $A$-module, $B \ot_A\Omt_A$ is a free rank $n+1$ 
$B$-module, and $\Omt_B \cong (B \ot_A\Omt_A) / \alpha(I / I^2)$.  For each $K$-valued
point $b \in X$, $(X\prl_b)^\circ \sub (\tau \A^n_b)^\circ$ is the linear subspace 
equal to $\cap_{f \in I}\Ker(\tau f(b))$, as desired.
\epr

\msk

The following corollary follows easily, by the standard embeddings of 
$X\prl$ and $TX$ in affine space.

\bco
Using notation from the previous proposition, the prolongation $X\prl$ 
is naturally isomorphic to the intersection of $CX$ with the hyperplane of 
$\A^{2n+1}$ defined by $\tau e -1=0$.  
Likewise, the tangent space $TX$ is naturally isomorphic to the 
intersection of $CX$ with the hyperplane of $\A^{2n+1}$ defined by $\tau e =0$.

With these natural embeddings of $X\prl$ and $TX$ in $CX$, $TX$ is a 
vector subbundle of $CX$, and $X\prl$ is an affine subbundle of $CX$.
Further, $X\prl$ is a $TX$-torsor under the action given by vector space addition
in $CX$.

Finally, $TX$ and $X\prl$ are (disjoint) principal divisors of $CX$.
\eco

\brm
One can easily show that these embeddings, in local coordinates, are the same as
the ones described above in terms of surjective $R$-algebra homomorphisms.
\erm

\bibliography{tdiff.bib}
\bibliographystyle{alpha}

\end{document}